\documentclass{article}
\usepackage{amssymb,amsmath,amsthm}
\usepackage{amsxtra}
\usepackage[mathscr]{eucal}
\usepackage{graphics}
\usepackage{epsfig}

\renewcommand{\le}{\leqslant}
\renewcommand{\ge}{\geqslant}

\makeatletter                   
\renewcommand\@biblabel[1]{#1.} 
\makeatother

\newtheorem{theorem}{Theorem}[section]
\newtheorem{lemma}[theorem]{Lemma}

\theoremstyle{definition}

\newcommand{\myref}[1]{Equation~\eqref{#1}}
\newcommand{\mybinom}[2]{{\displaystyle\binom{#1}{#2}}}

\begin{document}
\title{The Continuing Story of Zeta}
\author{Graham Everest, Christian R\" ottger and Tom Ward}
\maketitle

\addtocounter{section}{1}
\noindent{\bf\thesection. TAKING THE LOW ROAD.}
\setcounter{theorem}{0}
Riemann's Zeta Function~$\zeta(s)$ is defined for complex~$s=\sigma
+it$ with~$\Re(s)=\sigma >1$ by the formula
$$\zeta(s)=\sum_{n=1}^{\infty}\frac{1}{n^s}.$$
There are many ways to obtain the analytic continuation
of~$\zeta(s)$ to the left hand half-plane. The high road, Riemann's
own~\cite{riemann}, uses contour integration at an early stage, and
leads directly to the functional equation. Many authors
(\cite{MR0434929}, \cite{GandT}, \cite{MR568909}, \cite{patterson},
\cite{prachar},
\cite{titchmarsh}, and~\cite{witandwat}) use this method, or
variants of it, often at a more leisurely pace. Other methods are
known (Chapter~$2$ of~\cite{titchmarsh} lists seven) but a toll
seems inevitable on any route ending with the functional equation.

There are lower roads which give both the continuation to the whole
plane and the evaluation at non-positive integers but stop short of
proving the functional equation. If these are rigorous, yet quick
and simple, there must surely be a case for using them as well. The point
of this article to draw wider attention to these, often very scenic,
roads. In
his beautiful article~\cite[Sect.~7]{ayoub}, Ayoub comments upon
Euler's paper of~$1740$ in which he boldly evaluates divergent
series to obtain~$\zeta(-k)$ for integers~$k\ge 0$, thereby
predicting the functional equation. Recently,
Sondow~\cite{sondow} has noted one way in which Euler's argument can
be made rigorous. Simultaneously, Min\'a\v c \cite{minac} showed
how to evaluate $\zeta(-k)$ in an extremely simple and elegant
way, by integrating a polynomial on~$[0,1]$. More recently,
Murty and Reece \cite{murtyreece} have shown how the continuation
and evaluation of the Hurwitz zeta function can be obtained in
a simple down-to-earth way and this is applicable to~$\zeta(s)$
and many $L$-functions. The point of this note is to highlight
just how easily the continuation and evaluation of~$\zeta(s)$ can
be obtained. All that we say can be found in the articles cited. For
example, our work-horse (\ref{zetarelation}) is the truncation of
Landau's formula \cite[p. 274]{landau}.

\bigskip

\addtocounter{section}{1}
\noindent{\bf\thesection. A JOURNEY OF A
THOUSAND MILES...}\setcounter{theorem}{0}
Notice that for~$\sigma >1$,
\begin{equation*}\label{1/s-1}\int_1^\infty x^{-s}\; dx =
\frac{-1}{1-s}=\frac{1}{s-1},
\end{equation*}
which yields at once the continuation to the whole complex plane of
the function represented by the integral for~$\sigma>1$. Obviously
the continuation is analytic everywhere apart from a simple pole
at~$s=1$. For~$\sigma>1$,
\begin{eqnarray}\label{equation:lowroad}
   \frac{1}{s-1} &=& \int_1^\infty x^{-s}\; dx\nonumber
       = \sum_{n=1}^\infty \int_n^{n+1} x^{-s}\; dx\nonumber\\
       &=& \sum_{n=1}^\infty \int_0^1 (n+x)^{-s}\; dx
       = \sum_{n=1}^\infty \frac{1}{n^s}
               \int_0^1\left(1+\frac{x}{n}\right)^{-s}\;dx.\label{integralsumtosubstitute}
\end{eqnarray}
All the sums converge absolutely for~$\sigma>1$. In what follows we
assume that~$\sigma >1$ and that~$|s|$ is bounded by~$K$, a fixed
(although arbitrary) constant. Now begin the binomial expansion of
the bracketed term, noting that the higher binomial coefficients all
include a factor~$s$:
\begin{equation}\label{equation:guile}
\left(1+\frac{x}{n}\right)^{-s}=1-\frac{sx}{n}+ sE_1(s,x,n),
\end{equation}
where the function~$E_1$ satisfies
\begin{equation}\label{equation:errorbound1}
   |E_1(s,x,n)|\le \frac{C_1x^2}{n^2}\le \frac{C_1}{n^2},
\end{equation}
for all~$x\in[0,1]$ and all~$n\ge1$, with~$C_1=C_1(K)$ (since~$E_1$
is just the error term of a Taylor series in~$x/n$).
Substitute~\myref{equation:guile} into the
sum~\eqref{integralsumtosubstitute} and perform the integration with
respect to~$x$. We find that
\begin{equation}
   \label{A1}
   \frac{1}{s-1}=\zeta(s) -\frac{s}{2}\zeta(s+1) +sA_1(s),
\end{equation}
where~$A_1(s)$ is analytic for~$\sigma>-1$
by~\eqref{equation:errorbound1}. Thus \myref{A1} may be used to
extend~$\zeta(s)$ to the half-plane~$\sigma>0$. It even shows that
the extended function will be analytic there apart from a simple
pole at~$s=1$ with residue~$1$. In other words, \myref{A1} implies
that
\begin{equation}  \label{equation:knownresidue}
   \lim_{s\to 1} (s-1)\zeta(s) = 1.
\end{equation}
\myref{equation:knownresidue} can also be written~$\lim_{s\to 0}
s\zeta(s+1) = 1$. Using this fact, and letting~$s\to 0^+$ in
\myref{A1}, we obtain
 \[
 -1=\zeta(0)-\frac12,
 \]
which yields the known value~$\zeta(0)=-1/2$.

The preceding argument begins with the binomial
estimate~\eqref{equation:guile}, finds the analytic continuation of
the zeta function to the half-plane~$\sigma>0$ and
evaluates~$\zeta(0)$. What happens if more terms of the binomial
expansion are included? An additional term in the binomial expansion
gives
$$\left(1+\frac{x}{n}\right)^{-s}=1-\frac{sx}{n}+\frac{s(s+1)x^2}{2n^2}+
(s+1)E_2(s,x,n);
$$
notice that the higher binomial coefficients all include a
factor~$(s+1)$. Here,~$E_2$ is a function which satisfies
\begin{equation*}\label{equation:errorbound2}
   |E_2(s,x,n)|\le \frac{C_2x^3}{n^3}\le \frac{C_2}{n^3},
\end{equation*}
for all $x\in[0,1]$ and all $n$, where~$C_2=C_2(K)$. Substituting
this into~\eqref{integralsumtosubstitute} and integrating as before
gives
\begin{equation}\label{A2}
     \frac{1}{s-1}=\zeta(s) -\frac{s}{2}\zeta(s+1)+
\frac{s(s+1)}{6}\zeta(s+2) +(s+1)A_2(s),
\end{equation}
where~$A_2(s)$ is analytic for~$\sigma >-2$. Thus, \myref{A2} may be
used to continue~$\zeta(s)$ to the half-plane~$\sigma >-1$. As
before, letting~$s \rightarrow -1^+$ and using
\myref{equation:knownresidue}, we obtain
$$-\frac{1}{2}=\zeta(-1)+\frac{1}{2}\zeta(0)-\frac{1}{6}=
\zeta(-1)-\frac{1}{4}-\frac{1}{6}
$$
yielding the known value $\zeta(-1)=-1/12$.


\bigskip

\addtocounter{section}{1}
\noindent{\bf\thesection. GENERAL METHOD.}
This method can be repeated in order to continue~$\zeta(s)$ further
and further to the left of the complex plane. Moreover, it yields
the explicit evaluation at the non-positive integers in terms of the
\emph{Bernoulli numbers}. The sequence of \emph{Bernoulli numbers
$\left(B_n\right)$} is defined via the generating function
\begin{equation}\label{bernoulligen}
   \frac{x}{e^x-1}=\sum_{n=0}^\infty B_n\frac{x^n}{n!}
\end{equation}
from which it is clear that all the~$B_n$ are rational numbers. We
need two well-known properties of this fascinating sequence which
are stated in the following lemma.

\begin{lemma}
With~$B_n$ defined by~\eqref{bernoulligen},
\begin{equation}
   \label{bernrelation}
\sum_{n=0}^{N-1}\mybinom{N}{n} B_n=0
\qquad\mbox{for all~$N>1$},
   \end{equation}
and
\begin{equation*}
   \label{vanishodd}
   B_n=0 \qquad\mbox{for all odd~$n\ge 3$}.
\end{equation*}
\end{lemma}

\begin{proof}
The relation (\ref{bernoulligen}) can be written
$$
(e^x-1)\sum_{n=0}^\infty B_n\frac{x^n}{n!}=x.
$$
For~$N>1$ the coefficient of~$x^N$ in the left-hand side is
\[
\sum_{m=0}^{N-1}\frac{1}{(N-m)!m!}B_m,
\]
which gives~\eqref{bernrelation} after multiplying by~$N!$. The
second statement follows from the fact that
\[
\frac{x}{e^x-1} + \frac{x}{2}  = \frac{x(1+e^x)}{e^x-1}
\]
is an even function.
\end{proof}

The recurrence relation (\ref{bernoulligen}) can be used to
calculate~$B_n$ inductively. The first few Bernoulli numbers are
given below.
\begin{equation*}
\begin{array}{c|ccccccccccc}
n&0&1&2&3&4&5&6&7&8&9&10\\
\hline
B_n&\vphantom{\sum^A}\hphantom{\medspace}1\hphantom{\medspace}&-\frac12&
\hphantom{\medspace}\frac16\hphantom{\medspace}&
\hphantom{\medspace}0\hphantom{\medspace}&-\frac{1}{30}&
\hphantom{\medspace}0\hphantom{\medspace}&
\hphantom{\medspace}\frac{1}{42}\hphantom{\medspace}&
\hphantom{\medspace}0\hphantom{\medspace}&-\frac{1}{30}&
\hphantom{\medspace}0\hphantom{\medspace}&
\hphantom{\medspace}\frac{5}{66}\hphantom{\medspace}
\end{array}
\end{equation*}

\begin{theorem}\label{main}
There is an analytic continuation of~$\zeta(s)$ to the entire
complex plane where it is analytic apart from a simple pole at~$s=1$
with residue~$1$. For all~$k\ge 1$,
\begin{equation}\label{equation:main}
   \zeta(-k)=-\frac{B_{k+1}}{k+1}.
\end{equation}
\end{theorem}

Note that~\myref{equation:main} is not true for~$k=0$ but our method
has already given us the special value~$\zeta(0)=-1/2$.

\begin{proof}[\sc Proof of Theorem \ref{main}]
The analytic continuation of
the zeta function to the half-plane~$\sigma
>-k$ arises in exactly the same way as before, by extracting an
appropriate number of terms of the binomial expansion and using
induction. For integral~$k\ge 0$ and~$\sigma>1$, this gives the
relation
\begin{eqnarray}
\frac{1}{s-1}&=&
   \zeta(s)+\sum_{r=0}^k \frac{(-1)^{r+1} s(s+1)\dots
(s+r)}{(r+2)!}\zeta(s+r+1)\nonumber\\
&&\qquad\qquad\qquad\qquad\qquad\qquad
+(s+k)A_{k+1}(s)\label{zetarelation}
\end{eqnarray}
where~$A_{k+1}(s)$ is analytic in~$\sigma>-(k+1)$, again because all
higher binomial coefficients include a factor~$(s+k)$. Notice
that~$k=0$ gives \myref{A1} and~$k=1$ gives \myref{A2}.

By induction, we may assume that $\zeta(s)$ has already been
extended to the half-plane $\sigma>1-k$ so \myref{zetarelation} is
valid there, because the singularities at $s=0,-1,\dots$ are
removable. All the functions in \myref{zetarelation} except
$\zeta(s)$ are defined at least for $\sigma>-k$, which gives the
analytic continuation of $\zeta(s)$ to that half-plane. Let
$s\rightarrow -k^+$ in~\eqref{zetarelation} and use
\myref{equation:knownresidue} for the term with~$r=k$ to obtain
$$
-\frac{1}{k+1}=\zeta(-k)+\sum_{r=0}^{k-1} \mybinom{k}{r+1}
\frac{\zeta(-k+r+1)}{r+2}- \frac{1}{(k+1)(k+2)}.
$$
Writing $r$ for every $r+1$ simplifies this to
$$0=\zeta(-k)+\frac{1}{k+2}+\sum_{r=1}^{k} \mybinom{k}{r}
\frac{\zeta(-k+r)}{r+1}.
$$
The term with $r=k$ is known. Using induction on the
others gives
\begin{equation}\label{equation:recurrence1}
   0=\zeta(-k)+\frac{1}{k+2}
           -\sum_{r=1}^{k-1} \mybinom{k}{r}
\frac{B_{k-r+1}}{(r+1)(k-r+1)}-\frac{1}{2(k+1)}.
\end{equation}
A simple manipulation of factorials gives
\begin{equation*}\label{equation:binomials}
\frac{(k+1)(k+2)}{(r+1)(k-r+1)} \mybinom{k}{r} = \mybinom{k+2}{r+1}
 = \mybinom{k+2}{k-r+1},
\end{equation*}
which transforms \myref{equation:recurrence1} to
\begin{equation}\label{equation:recurrence2}
   0=\zeta(-k)+\frac{k}{2(k+1)(k+2)}
           -\frac{1}{(k+1)(k+2)}\sum_{r=1}^{k-1} \mybinom{k+2}{k-r+1}B_{k-r+1}.
\end{equation}
Multiply by $(k+1)(k+2)$ and apply \myref{bernoulligen} with
$N=k+2$. Only the terms for $r=0,k,k+1$, missing in
\myref{equation:recurrence2} survive, yielding
\begin{eqnarray}\label{equation:recurrence3}
   0&=& (k+1)(k+2)\zeta(-k)+\frac{k}{2}
           + (k+2) B_{k+1} + (k+2) B_1 + B_0 \nonumber\\
           &=&  (k+1)(k+2)\zeta(-k) +  (k+2)B_{k+1}\nonumber
\end{eqnarray}
and this completes the induction argument.
\end{proof}

\smallskip

\noindent{\bf ACKNOWLEDGEMENTS.} Our thanks go to J\'an Min\'a\v c
and Robin Kronenberg for helpful comments.

\end{document}